\documentclass{llncs}
\usepackage{amsfonts}
\usepackage{amssymb}
\usepackage{graphicx}
\usepackage{psfrag}
\usepackage[all]{xy}

\begin{document}

\title{Persistent topology for natural data analysis --- A survey}

\titlerunning{Persistence for natural data}  
%
\author{Massimo Ferri}
\authorrunning{M. Ferri} 
%
%
\institute{Dip. di Matematica and ARCES, Univ. di Bologna, Italy\\
\email{massimo.ferri@unibo.it},
}

\maketitle              

\begin{abstract}

Natural data offer a hard challenge to data analysis. One set of tools is being developed by several teams to face this difficult task: Persistent topology. After a brief introduction to this theory, some applications to the analysis and classification of cells, liver and skin lesions, music pieces, gait, oil and gas reservoirs, cyclones, galaxies, bones, brain connections, languages, handwritten and gestured letters are shown.
\end{abstract}


\section{Introduction and motivation}

What is the particular challenge offered by natural data, which could suggest the need of topology, and in particular of persistence? Simply said, it's quality instead of quantity. This is especially evident with images.

If one has to analyze, classify, retrieve images of mechanical pieces, vehicles, rigid objects, then geometry fulfills all needs. On the images themselves, matrix theory provides the transformations for superimposing a picture to a template. More often, pictures are represented by feature vectors, whose components are geometric measures ({\it shape descriptors}). Then recognition, defect detection, retrieval etc. can be performed on the feature vectors.

The scene changes if the depicted objects are of natural origin: the rigidity of geometry becomes an obstacle. Recognizing the resemblance between a sitting and a standing man is difficult. The challenge is even harder when it comes to biomedical data and when the context is essential for the understanding of data \cite{Ho12,ZiHiHo12}.

It's here that topology comes into play: the standing and sitting men are {\it homeomorphic}, i.e. there is a topological transformation which superimposes one to the other, whereas no matrix will ever be able to do that. It is generally difficult to discover whether two objects are homeomorphic; then algebraic topology turns helpful: It associates invariants --- e.g. Betti numbers --- to topological spaces, such that  objects which are homeomorphic have identical invariants (the converse does not hold, unfortunately).

(Algebraic) topology seems then to be the right environment for formalizing qualitative aspects in a computable way, as is nicely expressed in \cite[Sect. 5.1]{Ho14}. There is a problem: if geometry is too rigid, topology is too free. This is the reason why persistent topology can offer new topological descriptors (e.g. Persistent Betti Numbers, Persistence Diagrams) which preserve some selected geometric features through {\it filtering functions}. Classical references on persistence are \cite{VeUrFrFe93,CaZo*05,BiCeFrGiLa08,EdHa08}.

Persistent topology has been experimented in the image context, particularly in the biomedical domain, but also in fields where data are not pictures, e.g. in geology,  music and linguistics, as will be shown in this survey.

\section{Glossary and basic notions}

It is out of the scope of this survey to give a working introduction to homology and persistence; we limit ourselves to an intuitive description of the concepts, and recommend to profit of the technical references, without which a real understanding of the results is impossible. An essential (and avoidable) technical description of a particular homology is reported in Section \ref{tech}.

\medskip
\noindent {\bf Homology} There is a well-structured way (technically a set of functors) to associate {\it homology vector spaces} (more generally modules) $H_k(X)$ to a simplicial complex or to a topological space $X$, and linear transformations to maps \cite[Ch. 2]{Ha01} \cite[Ch. 4]{EdHa09}.

\medskip
\noindent {\bf Betti numbers} The {\it $k$-th Betti number} $\beta_k(X)$ is the dimension of the $k$-th homology vector space $H_k(X)$, i.e. the number of independent generators (homology classes of  {\it $k$-cycles}) of this space. Intuitively, $\beta_0(X)$ counts the number of path-connected components (i.e. the separate pieces) of which $X$ is composed; $\beta_1(X)$ counts the holes of the type of a circle (like the one of a doughnut); $\beta_2(X)$ counts the 2-dimensional voids (like the ones of gruyere or of an air chamber).

\medskip
\noindent {\bf Homeomorphism} Given topological spaces $X$ and $Y$, a {\it homeomorphism} from $X$ to $Y$ is a continuous map with continuous inverse. If one exists, the two spaces are said to be {\it homeomorphic}. This is the typical equivalence relation between topological spaces. Homology vector spaces and Betti numbers are invariant under homeomorphisms.

\begin{figure}[htb]
\begin{center}
  \includegraphics[width=1.0\textwidth]{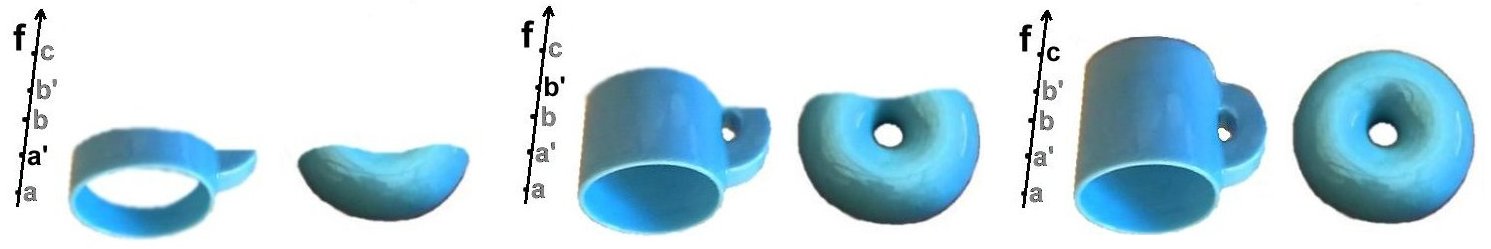}
    \caption{Sublevel sets of mug and doughnut.}
    \label{mug}
    \end{center}
\end{figure}

\begin{remark}
As hinted in the Introduction, geometry is too rigid, but topology is too free. In particular, homeomorphic spaces can be very different from an intuitive viewpoint: the joke by which ``for a topologist a mug and a doughnut are the same'' is actually true; the two objects are homeomorphic! {\bf Persistent topology} then tries to overcome this difficulty by studying not just topological spaces but pairs, once called {\it size pairs}, $(X, f)$ where $f$ is generally a continuous function, called {\it measuring} or {\it filtering function}, from $X$ to $\mathbb{R}$ (to $\mathbb{R}^n$ in {\it multidimensional persistence}) which conveys the idea of shape, the viewpoint of the observer. Shape similarity is actually very much dependent on the context. The Betti numbers of the sublevel sets then make it possible to distinguish the two objects although they are homeomorphic: see Figure \ref{mug}.
\end{remark}

\noindent {\bf Sublevel sets} Given a pair $(X, f)$, with $f:X \to \mathbb{R}$ continuous, given $u \in \mathbb{R}$, the {\it sublevel set under $u$} is the set $X_u = \{x\in X \,| \, f(x) \le u\}$. 

\begin{figure}[htb]
\begin{center}
  \includegraphics[width=1.0\textwidth]{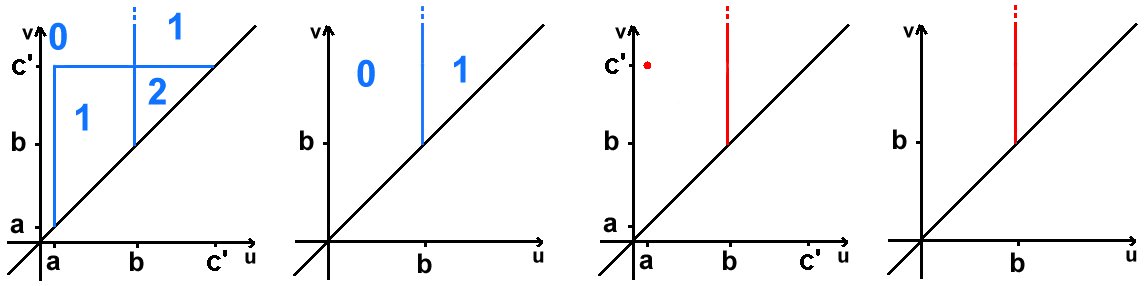}
    \caption{From left to right: 1-PBN functions of mug and of doughnut, 1-PDs of mug and of doughnut.}
    \label{betti}
    \end{center}
\end{figure}

\medskip
\noindent {\bf Persistent Betti Numbers} For all $u, v \in \mathbb{R},\ \ \ u<v$, the inclusion map $\iota^{u, v}: X_u \to X_v$ is continuous and induces, at each degree $k$, a linear transformation $\iota^{u, v}_*: H_k(X_u) \to H_k(X_v)$. The {\it $k$-Persistent Betti Number ($k$-PBN) function} assigns to the pair $(u, v)$ the number dim Im $\iota^{u, v}_*$, i.e. the number of classes of $k$-cycles of $H_k(X_u)$ which ``survive'' in $H_k(X_v)$. See Figure \ref{betti} (left) for the 1-PBN functions of mug and doughnut. Note that a pitcher, and more generally any open container with a handle, will have very similar PBNs to the ones of the mug; this is precisely what we want for a functional search and not for a strictly geometrical one.

\noindent {\bf Persistence Diagrams}. The $k$-PBN functions are wholly determined by the position of some discontinuity points and lines, called {\it cornerpoints} and {\it cornerlines} (or {\it cornerpoints at infinity}) The coordinates $(u, v)$ of a cornerpoint represent the levels of ``birth'' and ``death'' respectively of a generator; the abscissa of a cornerline is the level of birth of a generator which never dies. The {\it persistence} of a cornerpoint is the difference $v-u$ of its coordinates. Cornerpoints and cornerlines form the {\it $k$-Persistence Diagram} ($k$-{\it PD}). Figure \ref{betti} (right) depicts the  1-PDs of mug and of doughnut.
For the sake of simplicity, we are here neglecting the fact that cornerpoints and cornerlines may have multiplicities.

\begin{figure}[htb]
\begin{center}
  \includegraphics[width=1.0\textwidth]{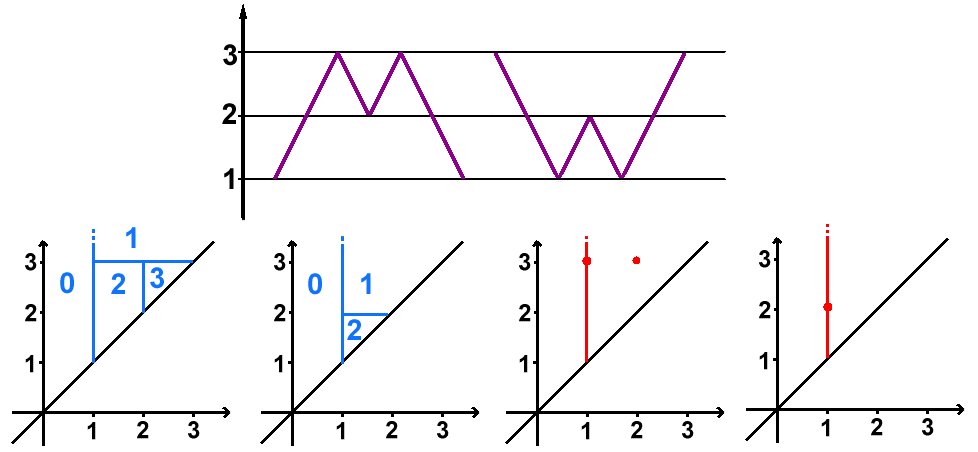}
    \caption{Above: the objects ``M'' and ``W''. Below, from left to right: 0-PBN functions of M and of W, 0-PDs of M and of W.}
    \label{MW}
    \end{center}
\end{figure}

\begin{remark}
Sometimes it is important to distinguish even objects for which there exists a rigid movement superimposing one to the other --- so also geometrically equivalent --- as in the case of some letters: context may be essential! See Figure \ref{MW}, where ordinate plays the role of filtering function.
\end{remark}

\noindent {\bf Matching distance} Given the $k$-PDs ${\cal D}_{X,f}, {\cal D}_{Y,g}$ of two pairs $(X, f), (Y, g)$, match the cornerpoints of ${\cal D}_{X,f}$ either with cornerpoints of ${\cal D}_{Y,g}$ or with their own projections on the diagonal $u=v$; the {\it weight} of this matching is the sup of the $L_\infty$-distances of matching points. The {\it matching distance} (or {\it bottleneck distance}) of ${\cal D}_{X,f}$ and ${\cal D}_{Y,g}$ is the inf of such weights among all possible such matchings.

\medskip
\noindent {\bf Natural pseudodistance} Given two pairs $(X, f), (Y, g)$, with $X, Y$ homeomorphic, the {\it weight} of a given homeomorphism $\varphi:X \to Y$ is sup$_{x\in X}|g(\varphi(x))-f(x)|$. The {\it natural pseudodistance} of $(X, f)$ and $(Y, g)$ is the inf of these weights among all possible homeomorphisms. If we are given the $k$-PDs of the two pairs, their  matching distance is a lower bound for the natural pseudodistance of the two pairs, and it is the best possible obtainable from the two $k$-PDs. Much is known on this dissimilarity measure \cite{DoFr04,DoFr04bis,DoFr07}.

\subsection{A brief technical description of homology}\label{tech}

There are several homologies. The classical and most descriptive one, at least for compact spaces, is singular homology with coefficients in $\mathbb{Z}$; we refer to \cite[Ch. 2]{Ha01} for a thorough exposition of it. Anyway, the homology used in most applications is the simplicial one, of which (with coefficients in $\mathbb{Z}_2$) we now give a very short introduction following \cite[Ch. 4]{EdHa09}.

\medskip
\noindent {\bf Simplices} A {\it $p$-simplex} $\sigma$ is the convex hull, in a Euclidean space, of a set of $p+1$ points, called {\it vertices} of the simplex, not contained in a Euclidean $(p-1)$-dimensional subspace; the simplex is said to be {\it generated} by its vertices. A {\it face} of a simplex $\sigma$ is the simplex generated by a nonempty set of vertices of $\sigma$.

\medskip
\noindent {\bf Simplicial complexes} A finite collection $K$ of simplices of a given Euclidean space is a {\it simplicial complex} if 1) for any $\sigma \in K$, all faces of $\sigma$ belong to $K$, 2) the intersection of two simplices of $K$ is either empty or a common face. The {\it space} of the complex $K$ is the topological subspace of Euclidean space $|K|$ formed by the union of all simplices of $K$.

\medskip
\noindent {\bf Simplicial homology with $\mathbb{Z}_2$ coefficients}. Given a (finite) simplicial complex $K$, call {\it $p$-chain} any formal linear combination of $p$-simplices with coefficients in $\mathbb{Z}_2$ (i.e. either 1 or 0, with $1+1=0$). $p$-chains form a $\mathbb{Z}_2$-vector space $C_p$. Note that each $p$-chain actually identifies a set of $p$-simplices of $K$ and that the sum of two $p$-chains is just the symmetric difference (Xor) of the corresponding sets. We now introduce a linear transformation $\partial_p: C_p \to C_{p-1}$ (called {\it boundary operator}) for any $p\in \mathbb{Z}$. We just need to define it on generators, i.e. on $p$-simplices, and then extend by linearity. Writing $\sigma = [u_0, u_1, \ldots, u_p]$, we denote by $[u_0, \ldots, \hat{u}_j, \ldots, u_p]$ its face generated by all of its vertices except $u_j$ ($j=0, \ldots, p$). Then we define
$$\partial_p(\sigma) = \sum_{j=0}^n[u_0, \ldots, \hat{u}_j, \ldots, u_p]$$

\begin{figure}[htb]
\begin{center}
  \includegraphics[width=0.6\textwidth]{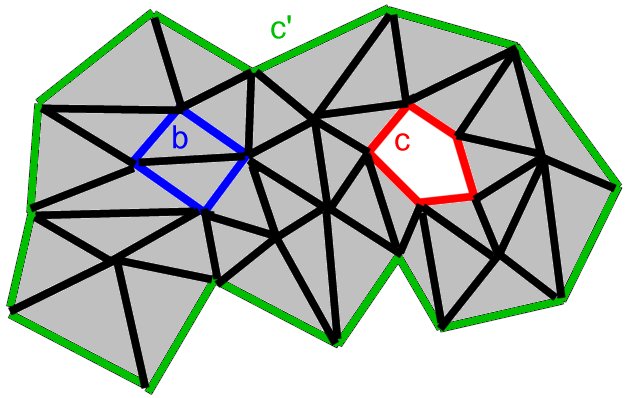}
    \caption{Cycles.}
    \label{cycle}
    \end{center}
\end{figure}

It is possible to prove that $\partial_{p}\partial_{p+1} = 0$, so that $B_p=$ Im$\partial_{p+1}$ is contained in $Z_p=$ Ker$\partial_p$. Elements of $B_p$ are called {\it $p$-boundaries}; elements of $Z_p$ are called {\it $p$-cycles}. The {\it $p$-homology vector space} is defined as the quotient  $H_p(K) = Z_p/B_p$. Homology classes are represented by cycles which are not boundaries. Two cycles are {\it homologous} is their difference is a boundary. In Figure \ref{cycle}, representing the simplicial complex $K$ formed by the shaded triangles and their faces, the blue chain $b$ is a 1-cycle which is also a boundary; the red chain $c$ and the green one $c'$ are 1-cycles which are not boundaries; $c$ and $c'$ are homologous.

\section{State-of-the-art}

The application of persistence to shape analysis and classification has a long story, since it started in the 90's when it still had the name of Size Theory \cite{VeUrFrFe93}. In the last few years it has taken various, very interesting forms. The constant aspect is always the presence of qualitative features which are difficult to capture and formalize within other frames of mind.

\subsection{Leukocytes}\label{sec:leuko}

Leukocytes, or white blood cells, belong to five different classes: lymphocyte; monocyte; neutrophile, eosinophile, basophile granulocytes. 
Eosino\-phile and neutrophile granulocytes are generally difficult to be distinguished, so they were considered in a single classification class in an early research by the Bologna team \cite{FeLoPa94}.

\begin{figure}[htb]
  \centering
  \begin{minipage}[ht]{0.3\linewidth}
     \centering
      \includegraphics[width=0.65\linewidth]{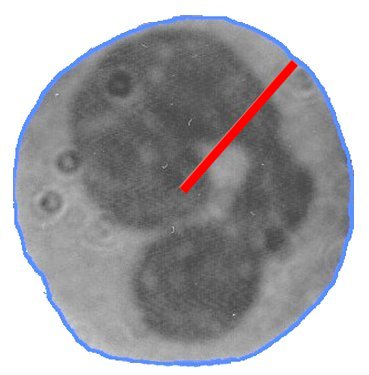}
    \caption{A radius along which the three filtering functions are computed.}
    \label{radius}
  \end{minipage}%
  \hfill%
  \begin{minipage}[ht]{0.7\linewidth}
     \centering
     \includegraphics[width=0.95\linewidth]{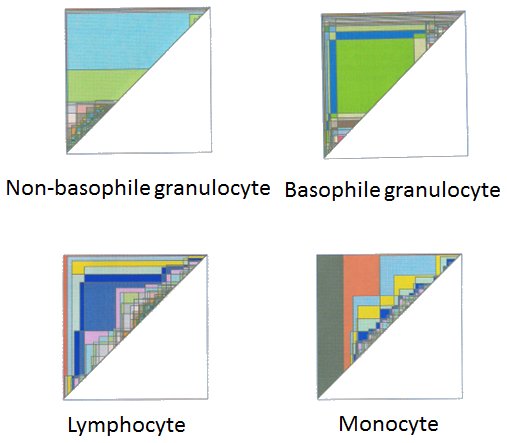}
     \caption{Persistent Betti Number functions relative to the sum of grey tones (different colors represent different values).} \label{leukobetti}
  \end{minipage}
\end{figure}

As a space, the boundary of the starlike hull of the cell is assumed. The images are converted to grey tones.

Three filtering functions are put to work, all computed along radii from the center of mass of the cell (Figure \ref{radius}):
\begin{itemize}
\item Sum of grey tones
\item Maximum variation
\item Sum of variations pixel to pixel.
\end{itemize}

Classification (with very good hit ratios for that time) is performed by measuring distance from the average PBN function of each class.

\subsection{Handwritten letters and monograms}\label{letters}

Again in Bologna we faced recognition of handwritten letters with time information; our goal was to recognize both the alphabet letter and the writer \cite{FeGa*95}.

The space on which the filtering functions are defined is the time interval of the writing. The filtering functions are computed in the 3D ``plane-time'':
\begin{itemize}
\item Distance of points from the letter axis
\item Speed
\item Curvature
\item Torsion
\item Distance from center of mass (in plane projection).
\end{itemize}

\begin{figure}[htb]
\begin{center}
  \includegraphics[width=0.8\textwidth]{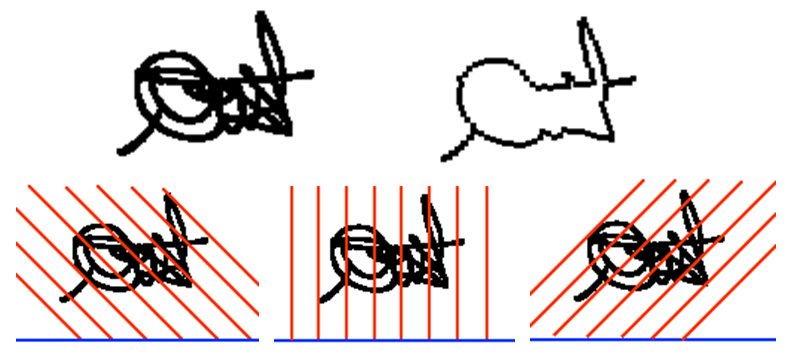}
    \caption{A monogram with its outline (above) and the directions along which the filtering functions are computed (below).}
    \label{monogram}
    \end{center}
\end{figure}

Classification comes from fuzzy characteristic functions, obtained from normalized inverse of distance. Cooperation of the characteristic functions coming from the single filtering functions is given by their rough arithmetic average.

\medskip
A later experiment, which was even repeated live at a conference, concerned the recognition of monograms for personal identification, without time information \cite{FeFr*98}.

Two topological spaces are used. The first is the outline of the monogram and the filtering function is the distance from the center of mass (see upper Figure \ref{monogram}).

The second space is a horizontal segment placed at the base of the monogram image. Filtering functions:

\begin{itemize}
\item Number of black pixels along segments (3 directions) (see lower Figure \ref{monogram})
\item Number of pixel-pixel black-white jumps (3 directions).
\end{itemize}

Classification is performed by a weighted average of fuzzy characteristic functions.

\subsection{Sign alphabet}\label{sign}

Automatic recognition of the symbols expressed by the hands in the sign language is a task which was of interest for different teams. The first one was the group led by Alessandro Verri in Genova \cite{UrVe94}. The signs were performed with a white glove on a black background;
translation into common letters was done in real time in a live demo at a conference.


The domain space is a horizontal segment; the filtering functions assign to each point of the segment the maximum distance of a contour point within a strip of fixed width, with 24 different strip orientations.

\begin{figure}[htb]
\begin{center}
  \includegraphics[width=1.0\textwidth]{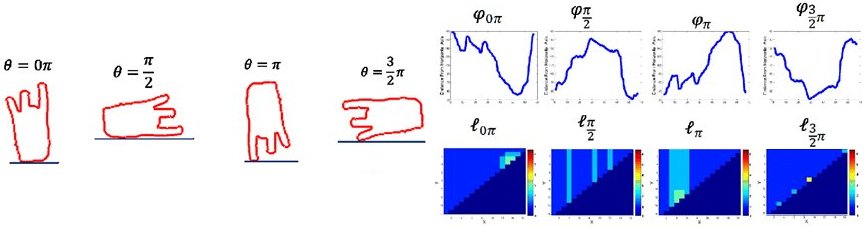}
    \caption{Four filtering functions and the corresponding 0-Persistent Betti Number functions.}
    \label{sign2}
    \end{center}
\end{figure}

The choice of S. Wang in Sherbrooke, instead, is to use a part of the contour, determined by principal component analysis, as a domain and distance from center of mass as filtering function \cite{HaZiWa99}.

\medskip
The team of D.Kelly in Maynooth uses the whole contour as domain, and distances from four lines as filtering functions \cite{KeMc*08} (see Figure \ref{sign2}).

\subsection{Human gait}\label{gait}

Personal identification and surveillance are the aim of a research by the Cuban team of L. Lamar-Le\'on, together with the Sevilla group of computational topology \cite{LaGa*12}.

\begin{figure}[htb]
\begin{center}
  \includegraphics[width=1.0\textwidth]{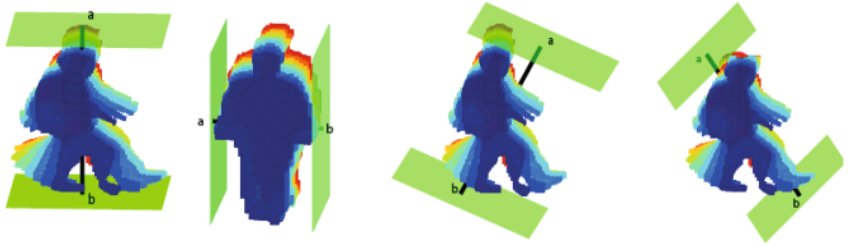}
    \caption{Four filtering functions on silhouette stacks for gait identification.}
    \label{fig:gait}
    \end{center}
\end{figure}

Considering a stack of silhouettes as a 3D object, and using four different filtering functions, makes 0- and 1-degree persistent homology a tool for identifying people through their gait (Figure \ref{fig:gait}).

\subsection{Tropical cyclones}\label{cycl}

\begin{figure}[htb]
\begin{center}
  \includegraphics[width=1.0\textwidth]{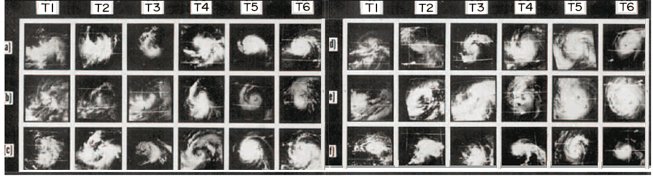}
    \caption{Time evolution of cyclones.}
    \label{fig:cycl}
    \end{center}
\end{figure}

S. Banerjee in Kolkata makes use of persistence on sequences of satellite images of cloud systems (Figure \ref{fig:cycl}), in order to evaluate risk and intensity of forming hurricanes \cite{Ba11}.

Time interval is the domain of two filtering functions which are common characteristic measures of cyclones:
\begin{itemize}
\item Central Feature portion
\item Outer Banding Feature
\end{itemize}

\subsection{Galaxies}\label{galax}

Again S. Banerjee \cite{Ba14} applies similar methods to another type of spirals: galaxies.


Various filtering functions are used. One is defined as a function of distance from galaxy center, and is the ratio between major and minor axis of the corresponding isophote.
Another one is a ``pitch'' parameter defined by Ringermacher and Mead \cite{RiMe09}.
A third filtering function is a compound based on color.

The classification results agree with the literature.

\subsection{Bones}\label{bones}

In \cite{TuMu*14} a powerful construction (the Persistent Homology Transform) is introduced. It consists in gathering the ``height'' filtering functions according to all possible directions. The paper shows that the transform is injective for objects homeomorphic to spheres. By using the transform it is possible to define an effective distance between surfaces. An application is shown by classifying heel bones of different species; the comparison with the ground truth produced by using placement of landmarks on the surfaces is very good.

\subsection{Melanocytic lesions}\label{mela}

A very important part of natural shape analysis is the detection of malignant cells and lesions, since there generally are no templates for them. As far as we know, the first attempt through persistence (called {\it size theory} at that time) is the ADAM EU Project, by the Bologna team together with CINECA and with I. Stanganelli, a dermatologist of the Romagna Oncology Institute \cite{dAFe*04,StBr*05,FeSt10}.
The analysis is mainly based on asymmetry of boundary, masses and color distribution: the lesion is split into two halves by 45 equally spaced lines, and the difference between the two halves is measured by the matching distance of the corresponding Persistence Diagrams.

\begin{figure}[htb]
\begin{center}
  \includegraphics[width=1.0\textwidth]{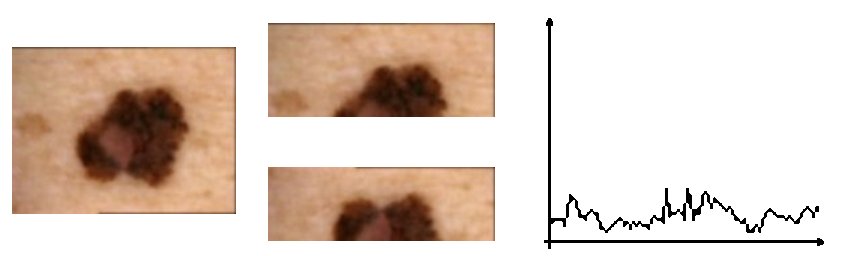}
    \caption{One of the 45 splittings of a melanocytic lesion, and the whole $A$-curve corresponding to the filtering function luminance.}
    \label{fig:mela}
    \end{center}
\end{figure}

The three functions ($A$-curves) relating these distances to the splitting line angles give parameters which are then fed into a Support Vector Machine classifier.

\medskip
The same team is presently involved with a biomedical firm in the realization of a machine for smart retrieval of dermatological images \cite{FeTo*16}.

\subsection{Tumor mouth cells}

\begin{figure}[htb]
\begin{center}
  \includegraphics[width=0.6\textwidth]{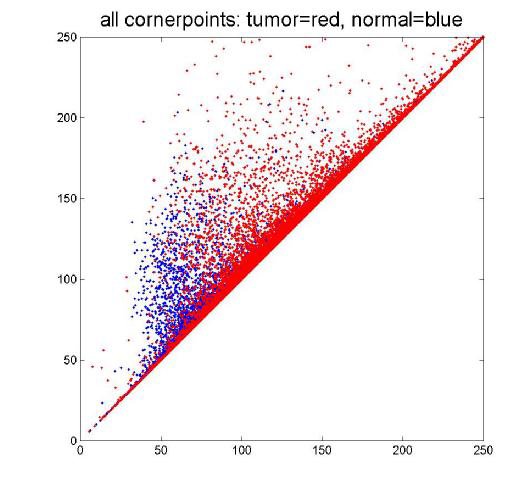}
    \caption{Distribution of cornerpoints in the diagrams of normal and tumor mouth cells.}
    \label{fig:mouth}
    \end{center}
\end{figure}

A morphological classification of normal and tumor cells of the epithelial tissue of the mouth is proposed in \cite{Mi06,MiLa08}: the filtering function is distance from the center of mass; the discrimination is statistically based on the distribution of cornerpoints (see Fig.~\ref{fig:mouth}).

\subsection{Hepatic lesions}\label{hepa}

\begin{figure}[htb]
\begin{center}
  \includegraphics[width=0.88\textwidth]{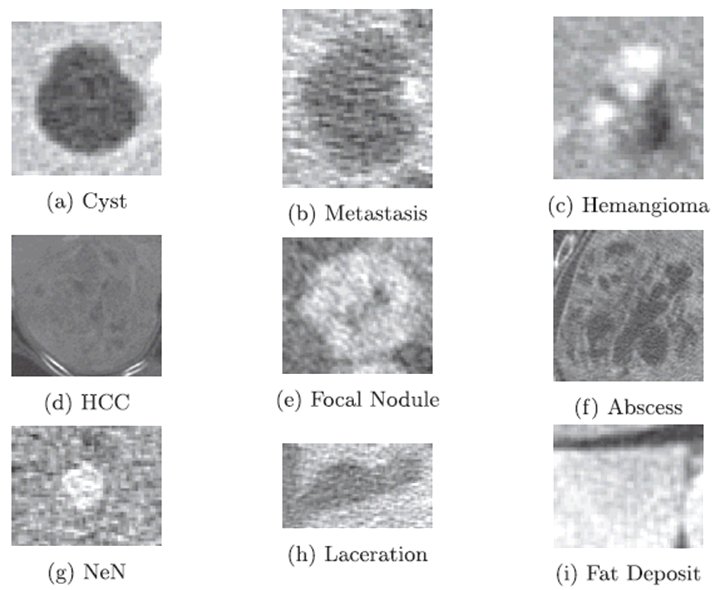}
    \caption{Various types of hepatic lesions.}
    \label{fig:hepa}
    \end{center}
\end{figure}

The advantages of a multidimensional range for the filtering functions are shown in \cite{AdRuCa14}, where several classification experiments are performed on the images of hepatic cells (see Figure \ref{fig:hepa}). The domain space is the part of image occupied by the lesion; the two components of the filtering function are the greyscale of each pixel and the distance from the lesion boundary.

\subsection{Genetic pathways}\label{genetic}

So far we have seen applications of persistence to  images of natural origin. But the modularity of the method opens the possibility to deal with data of very different nature. A first example is given by \cite{PlBa*16}, where persistence is used on the Vietoris-Rips complex in a space where points are {\em complex phenotypes} related together by the {\em Jaccard distance}. This made it possible to find systematic associations among metabolic syndrome variates that show distinctive genetic
association profiles.

\subsection{Oil and gas reservoirs}

Researchers in Ufa and Novosibirsk need to get a reliable geological and hydrodynamical model of gas and oil reservoirs out of noisy data; the model has to be robust under small perturbations. The authors have found an answer in persistent 0-, 1- and 2-cycles. The domain space is the 3D reservoir bed, and the filtering function is permeability, obtained as a decreasing function of radioactivity \cite{BaBa*13} (Russian; translated and completed in this same volume).

\subsection{Brain connections}\label{brain}

A complex research on brain connections and their modification under the assumption of a psychoactive substance (psilocybine) is performed in \cite{PeEx*14}  and extended in \cite{LoEx*16}. The construction starts with a complete graph whose vertices are cortical or subcortical regions; these, and their functional connectivity (expressed as weights on the edges) come from an elaborate processing of functional MRI data. Then the simplicial complex is built, whose simplices are the cliques (complete subgraphs) of the graph.

\begin{figure}[htb]
\begin{center}
  \includegraphics[width=0.9\textwidth]{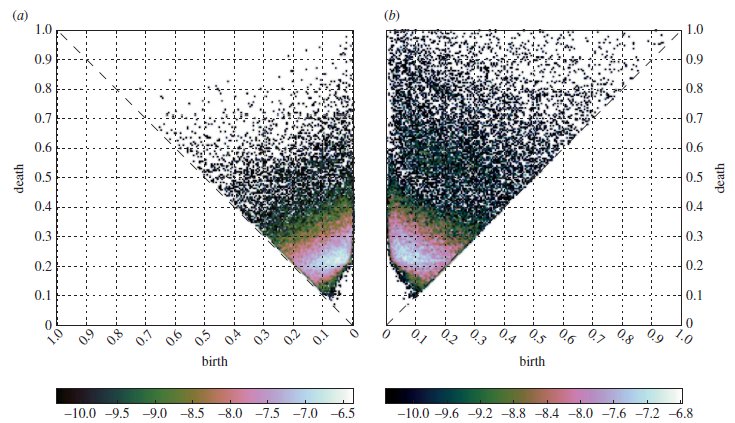}
    \caption{Probability densities for $H_1$ generators: placebo (left) and psilocybin (right) treated.}
    \label{fig:brain}
    \end{center}
\end{figure}

The filtering function on each simplex is minus the highest weight of its building edges. A difference between treated and control  subjects already appears in the comparison of the 1-Persistence Diagrams (see Figure \ref{fig:brain}). Then more information is obtained from secondary graphs (called {\it homological scaffolds}), whose vertices are the homology generators weighted by their persistence.

There are other applications of persistence to brain research: evaluation of cortical thickness in autism \cite{ChBu*09}; study of unexpected connections between subcortex, frontal cortex and parietal cortex in the form of 1- and 2-dimensional persistent cycles \cite{GiPa*15,SiGi*16}.

\subsection{Music}\label{music}

\begin{figure}[htb]
\begin{center}
  \includegraphics[width=0.9\textwidth]{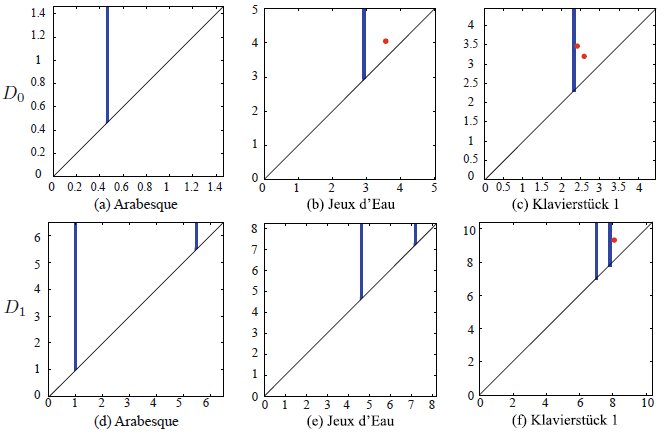}
    \caption{0- and 1-persistence diagrams for three classical pieces.}
    \label{fig:music}
    \end{center}
\end{figure}

Among other mathematical applications to music, M.G. Bergomi in Lisbon collaborates with various researchers in exploring musical genres by persistence \cite{BeBa*16}. As a space they adopt a modified version of Euler's {\it Tonnetz} \cite{BiAn*13}. The filtering function is the total duration of each note in a given track. Classification can be performed at different detail levels: experimentation is reported on tonal and atonal classical music of several authors (an example is in Figure \ref{fig:music}), on pop music and on different interpretation of the same jazz piece. 

A blend of persistence and deep learning is the central idea of a research by the team of I.-H. Yang in Taiwan \cite{LiJe*16}. They input audio signals to a Convolutional Neural Network (CNN); after a first convolution layer, a middle layer processes the output of the first in two different complementary ways: one is a classical CNN; the other computes the persistence landscape (an information piece derivable from the persistence diagram \cite{BuDl17}) of the same output. Whereas the persistence layer by itself does not perform any better than the normal CNN, their combination gives very good results in terms of music tagging.

\subsection{Languages}

An interdisciplinary team at Caltech investigates the metric spaces built by 79 Indo-European and 49 Niger-Congo languages \cite{PoGh*15}. These appear as points in a Euclidean space of syntactic parameters; on them a Vietoris-Rips complex \cite[Sect. III.2]{EdHa09} is built and Euclidean distance is assumed as filtering function. The Indo-European family reveals one 1-dimensional and two 0-dimensonal persistent cycles, the Niger-Congo respectively none and one. The interpretation of these differences and of the link with phylogenetic and historical facts is still under way.

\section{Open problems}

There is a number of open problems in persistence, whose solution will affect applications to natural data analysis, and to which only partial answers have been given so far:

\begin{itemize}
\item Optimal choice of the foliations along which to perform the 1D reduction of multidimensional persistence \cite{CeDiFeFrLa09}
\item Study of the discontinuities in multidimensional persistence \cite{CaZo09,CeFr15}
\item Understanding the monodromy around multiple cornerpoints \cite{CeEt*13}
\item Restricting the group of homeomorphisms of interest by considering the invariance required by the observer \cite{FrJa16}
\item Modulation of the impact of different filtering functions for search engines with relevance feedback \cite{GiFrSpFa10}
\item Use of advanced tools of algebraic topology \cite{BeMu15}
\item Use of persistence in the wider context of concrete categories, not necessarily passing through homology of complexes or of topological spaces \cite{BeFe*17}.
\end{itemize}

\section{Future outlook}

There are at least two ways in which persistence will interact with machine learning, and this is likely to enormously boost the qualitative processing of natural data \cite{DeEm*16}:

\begin{itemize}
\item{} Feeding a neural network with Persistence Diagrams instead of raw data will convey the needs and viewpoints of the user
\item{} Deep learning might yield a quantum leap in persistence, by automatically finding the best filtering functions for a given problem.
\end{itemize}

\section*{Acknowledgments}
Article written within the activity of INdAM-GNSAGA.

\bibliographystyle{abbrv}
\bibliography{FerriBIRS}

\begin{thebibliography}{10}

\bibitem{AdRuCa14}
A.~Adcock, D.~Rubin, and G.~Carlsson.
\newblock Classification of hepatic lesions using the matching metric.
\newblock {\em Computer vision and image understanding}, 121:36--42, 2014.

\bibitem{Ba11}
S.~Banerjee.
\newblock Size functions in the study of the evolution of cyclones.
\newblock {\em International Journal of Meteorology}, 36(358):39, 2011.

\bibitem{Ba14}
S.~Banerjee.
\newblock Size functions in galaxy morphology classification.
\newblock {\em Int. J. Comput. Appl}, 100(3):1--4, 2014.

\bibitem{BaBa*13}
Y.~V. Bazaikin, V.~A. Baikov, I.~A. Taimanov, and A.~A. Yakovlev.
\newblock Chislennyi analiz topologicheskih harakteristik trehmernyh
  geologicheskih modelei neftegazovyh mestorozhdenii.
\newblock {\em Matematicheskoe modelirovanie}, 25(10):19--31, 2013.

\bibitem{BeMu15}
F.~Belch{\'\i} and A.~Murillo.
\newblock A$_\infty$-persistence.
\newblock {\em Applicable Algebra in Engineering, Communication and Computing},
  26(1-2):121--139, 2015.

\bibitem{BeBa*16}
M.~G. Bergomi, A.~Barat{\`e}, and B.~{D}i {F}abio.
\newblock Towards a topological fingerprint of music.
\newblock In {\em International Workshop on Computational Topology in Image
  Context}, pages 88--100. Springer, 2016.

\bibitem{BeFe*17}
M.~G. Bergomi, M.~Ferri, and L.~Zuffi.
\newblock Graph persistence.
\newblock {\em arXiv preprint arXiv:1707.09670}, 2017.

\bibitem{BiCeFrGiLa08}
S.~Biasotti, A.~Cerri, P.~Frosini, D.~Giorgi, and C.~Landi.
\newblock Multidimensional size functions for shape comparison.
\newblock {\em Journal of Mathematical Imaging and Vision}, 32(2):161--179,
  2008.

\bibitem{BiAn*13}
L.~Bigo, M.~Andreatta, J.-L. Giavitto, O.~Michel, and A.~Spicher.
\newblock Computation and visualization of musical structures in chord-based
  simplicial complexes.
\newblock In {\em International Conference on Mathematics and Computation in
  Music}, pages 38--51. Springer, 2013.

\bibitem{BuDl17}
P.~Bubenik and P.~D{\l}otko.
\newblock A persistence landscapes toolbox for topological statistics.
\newblock {\em Journal of Symbolic Computation}, 78:91--114, 2017.

\bibitem{CaZo09}
G.~Carlsson and A.~Zomorodian.
\newblock The theory of multidimensional persistence.
\newblock {\em Discr. Comput. Geom.}, 42(1):71--93, 2009.

\bibitem{CaZo*05}
G.~Carlsson, A.~Zomorodian, A.~Collins, and L.~J. Guibas.
\newblock Persistence barcodes for shapes.
\newblock {\em IJSM}, 11(2):149--187, 2005.

\bibitem{CeDiFeFrLa09}
A.~Cerri, B.~{Di Fabio}, M.~Ferri, P.~Frosini, and C.~Landi.
\newblock Betti numbers in multidimensional persistent homology are stable
  functions.
\newblock {\em Mathematical Methods in the Applied Sciences},
  36(12):1543--1557, 2013.

\bibitem{CeEt*13}
A.~Cerri, M.~Ethier, and P.~Frosini.
\newblock A study of monodromy in the computation of multidimensional
  persistence.
\newblock In {\em International Conference on Discrete Geometry for Computer
  Imagery}, pages 192--202. Springer, 2013.

\bibitem{CeFr15}
A.~Cerri and P.~Frosini.
\newblock Necessary conditions for discontinuities of multidimensional
  persistent betti numbers.
\newblock {\em Mathematical Methods in the Applied Sciences}, 38(4):617--629,
  2015.

\bibitem{ChBu*09}
M.~Chung, P.~Bubenik, and P.~Kim.
\newblock Persistence diagrams of cortical surface data.
\newblock In {\em Information processing in medical imaging}, pages 386--397.
  Springer, 2009.

\bibitem{dAFe*04}
M.~d'Amico, M.~Ferri, and I.~Stanganelli.
\newblock Qualitative asymmetry measure for melanoma detection.
\newblock In {\em Biomedical Imaging: Nano to Macro, 2004. IEEE International
  Symposium on}, pages 1155--1158. IEEE, 2004.

\bibitem{DeEm*16}
M.~Dehmer, F.~Emmert-Streib, S.~Pickl, and A.~Holzinger.
\newblock {\em Big Data of Complex Networks}.
\newblock CRC Press, 2016.

\bibitem{DoFr04bis}
P.~Donatini and P.~Frosini.
\newblock Lower bounds for natural pseudodistances via size functions.
\newblock {\em Archives of {I}nequalities and {A}pplications}, 1(2):1--12,
  2004.

\bibitem{DoFr04}
P.~Donatini and P.~Frosini.
\newblock Natural pseudodistances between closed manifolds.
\newblock {\em Forum Mathematicum}, 16(5):695--715, 2004.

\bibitem{DoFr07}
P.~Donatini and P.~Frosini.
\newblock Natural pseudodistances between closed surfaces.
\newblock {\em Journal of the European Mathematical Society}, 9(2):231--253,
  2007.

\bibitem{EdHa08}
H.~Edelsbrunner and J.~Harer.
\newblock Persistent homology---a survey.
\newblock In {\em Surveys on discrete and computational geometry}, volume 453
  of {\em Contemp. Math.}, pages 257--282. Amer. Math. Soc., Providence, RI,
  2008.

\bibitem{EdHa09}
H.~Edelsbrunner and J.~Harer.
\newblock {\em Computational Topology: An Introduction}.
\newblock American Mathematical Society, 2009.

\bibitem{FeFr*98}
M.~Ferri, P.~Frosini, A.~Lovato, and C.~Zambelli.
\newblock Point selection: A new comparison scheme for size functions (with an
  application to monogram recognition).
\newblock In {\em Asian Conference on Computer Vision}, pages 329--337.
  Springer, 1998.

\bibitem{FeGa*95}
M.~Ferri, S.~Gallina, E.~Porcellini, and M.~Serena.
\newblock On-line character and writer recognition by size functions and fuzzy
  logic.
\newblock {\em Proc. ACCV’95}, pages 5--8, 1995.

\bibitem{FeLoPa94}
M.~Ferri, S.~Lombardini, and C.~Pallotti.
\newblock Leukocyte classifications by size functions.
\newblock In {\em Applications of Computer Vision, 1994., Proceedings of the
  Second IEEE Workshop on}, pages 223--229. IEEE, 1994.

\bibitem{FeSt10}
M.~Ferri and I.~Stanganelli.
\newblock Size functions for the morphological analysis of melanocytic lesions.
\newblock {\em Journal of Biomedical Imaging}, 2010:5, 2010.

\bibitem{FeTo*16}
M.~Ferri, I.~Tomba, A.~Visotti, and I.~Stanganelli.
\newblock A feasibility study for a persistent homology-based k-nearest
  neighbor search algorithm in melanoma detection.
\newblock {\em Journal of Mathematical Imaging and Vision}, pages 1--16, 2016.

\bibitem{FrJa16}
P.~Frosini and G.~Jab{\l}o{\'n}ski.
\newblock Combining persistent homology and invariance groups for shape
  comparison.
\newblock {\em Discrete \& Computational Geometry}, 55(2):373--409, 2016.

\bibitem{GiFrSpFa10}
D.~Giorgi, P.~Frosini, M.~Spagnuolo, and B.~Falcidieno.
\newblock 3{D} relevance feedback via multilevel relevance judgements.
\newblock {\em The Visual Computer}, 26(10):1321--1338, 2010.

\bibitem{GiPa*15}
C.~Giusti, E.~Pastalkova, C.~Curto, and V.~Itskov.
\newblock Clique topology reveals intrinsic geometric structure in neural
  correlations.
\newblock {\em Proceedings of the National Academy of Sciences},
  112(44):13455--13460, 2015.

\bibitem{HaZiWa99}
M.~Handouyahia, D.~Ziou, and S.~Wang.
\newblock Sign language recognition using moment-based size functions.
\newblock In {\em Proc. Int’l Conf. Vision Interface}, pages 210--216, 1999.

\bibitem{Ha01}
A.~Hatcher.
\newblock {\em Algebraic topology}.
\newblock Cambridge University Press, 2001.

\bibitem{Ho12}
A.~Holzinger.
\newblock On knowledge discovery and interactive intelligent visualization of
  biomedical data.
\newblock In {\em Proceedings of the Int. Conf. on Data Technologies and
  Applications DATA 2012, Rome (Italy)}, pages 5--16, 2012.

\bibitem{Ho14}
A.~Holzinger.
\newblock On topological data mining.
\newblock In {\em Interactive Knowledge Discovery and Data Mining in Biomedical
  Informatics}, pages 331--356. Springer, 2014.

\bibitem{KeMc*08}
D.~Kelly, J.~McDonald, T.~Lysaght, and C.~Markham.
\newblock Analysis of sign language gestures using size functions and principal
  component analysis.
\newblock In {\em Machine Vision and Image Processing Conference, 2008.
  IMVIP'08. International}, pages 31--36. IEEE, 2008.

\bibitem{LaGa*12}
J.~Lamar-Le{\'o}n, E.~B. Garcia-Reyes, and R.~Gonzalez-Diaz.
\newblock Human gait identification using persistent homology.
\newblock In {\em Iberoamerican Congress on Pattern Recognition}, pages
  244--251. Springer, 2012.

\bibitem{LiJe*16}
J.-Y. Liu, S.-K. Jeng, and Y.-H. Yang.
\newblock Applying topological persistence in convolutional neural network for
  music audio signals.
\newblock {\em arXiv preprint arXiv:1608.07373}, 2016.

\bibitem{LoEx*16}
L.-D. Lord, P.~Expert, H.~M. Fernandes, G.~Petri, T.~J. Van~Hartevelt,
  F.~Vaccarino, G.~Deco, F.~Turkheimer, and M.~L. Kringelbach.
\newblock Insights into brain architectures from the homological scaffolds of
  functional connectivity networks.
\newblock {\em Frontiers in Systems Neuroscience}, 10, 2016.

\bibitem{Mi06}
A.~Micheletti.
\newblock The theory of size functions applied to problems of statistical shape
  analysis.
\newblock In {\em S4G-International Conference in Stereology, Spatial
  Statistics and Stochastic Geometry}, pages 177--183. Union of Czech
  Mathematicians and Physicists, 2006.

\bibitem{MiLa08}
A.~Micheletti and G.~Landini.
\newblock Size functions applied to the statistical shape analysis and
  classification of tumor cells.
\newblock In {\em Progress in Industrial Mathematics at ECMI 2006}, pages
  538--542. Springer, 2008.

\bibitem{PeEx*14}
G.~Petri, P.~Expert, F.~Turkheimer, R.~Carhart-Harris, D.~Nutt, P.~J. Hellyer,
  and F.~Vaccarino.
\newblock Homological scaffolds of brain functional networks.
\newblock {\em Journal of The Royal Society Interface}, 11(101):20140873, 2014.

\bibitem{PlBa*16}
D.~E. Platt, S.~Basu, P.~A. Zalloua, and L.~Parida.
\newblock Characterizing redescriptions using persistent homology to isolate
  genetic pathways contributing to pathogenesis.
\newblock {\em BMC systems biology}, 10(1):S10, 2016.

\bibitem{PoGh*15}
A.~Port, I.~Gheorghita, D.~Guth, J.~M. Clark, C.~Liang, S.~Dasu, and
  M.~Marcolli.
\newblock Persistent topology of syntax.
\newblock {\em arXiv preprint arXiv:1507.05134}, 2015.

\bibitem{RiMe09}
H.~I. Ringermacher and L.~R. Mead.
\newblock A new formula describing the scaffold structure of spiral galaxies.
\newblock {\em Monthly Notices of the Royal Astronomical Society},
  397(1):164--171, 2009.

\bibitem{SiGi*16}
A.~Sizemore, C.~Giusti, R.~F. Betzel, and D.~S. Bassett.
\newblock Closures and cavities in the human connectome.
\newblock {\em arXiv preprint arXiv:1608.03520}, 2016.

\bibitem{StBr*05}
I.~Stanganelli, A.~Brucale, L.~Calori, R.~Gori, A.~Lovato, S.~Magi, B.~Kopf,
  R.~Bacchilega, V.~Rapisarda, A.~Testori, P.~A. Ascierto, E.~Simeone, and
  M.~Ferri.
\newblock Computer-aided diagnosis of melanocytic lesions.
\newblock {\em Anticancer Research}, 25(6C):4577--4582, 2005.

\bibitem{TuMu*14}
K.~Turner, S.~Mukherjee, and D.~M. Boyer.
\newblock Persistent homology transform for modeling shapes and surfaces.
\newblock {\em Information and Inference: A Journal of the IMA}, 3(4):310--344,
  2014.

\bibitem{UrVe94}
C.~Uras and A.~Verri.
\newblock On the recognition of the alphabet of the sign language through size
  functions.
\newblock In {\em Pattern Recognition, 1994. Vol. 2-Conference B: Computer
  Vision \& Image Processing., Proceedings of the 12th IAPR International.
  Conference on}, volume~2, pages 334--338. IEEE, 1994.

\bibitem{VeUrFrFe93}
A.~Verri, C.~Uras, P.~Frosini, and M.~Ferri.
\newblock On the use of size functions for shape analysis.
\newblock {\em Biol. Cybern.}, 70:99--107, 1993.

\bibitem{ZiHiHo12}
M.~Ziefle, S.~Himmel, and A.~Holzinger.
\newblock How usage context shapes evaluation and adoption criteria in
  different technologies.
\newblock In {\em AHFE 2012, Proceeding of Int. Conf. on Applied Human Factors
  and Ergonomics, San Francisco}, pages 2812--2821, 2012.

\end{thebibliography}

\end{document}